\documentclass[reqno,11pt]{amsart}

\newtheorem{theorem}{Theorem}[section]

\theoremstyle{definition}

\newtheorem{assumption}[theorem]{Assumption}

\theoremstyle{remark}
\newtheorem{remark}[theorem]{Remark}

\newcommand{\mysection}[1]{\section{#1}
\setcounter{equation}{0}}

\newcommand{\bR}{\mathbb R}

\newcommand\cF{\mathcal{F}}

\newcommand\cR{\mathcal{R}}

\newcommand{\HO}{\overset{\scriptscriptstyle\,0}%
 H \,\!}

\begin{document}
\title[Short proof
of It\^o's formula for   SPDEs]
{A relatively short proof
of It\^o's formula for  SPDEs and its applications}
 
\author{N.V. Krylov}
\address{127 Vincent Hall, University of Minnesota, Minneapolis,
 MN, 55455}
\thanks{The work   was partially supported by
NSF Grant  DNS-1160569}
\email{krylov@math.umn.edu}
 

\begin{abstract}
We give a short proof of It\^o's formula for stochastic
Hilbert-space valued processes in the setting $V\subset H\subset
V^{*}$ based on the possibility to lift the stochastic
differentials, which are originally in $V^{*}$,  into $H$.
Using this result we also prove the maximum
principle for second-order SPDEs in arbitrary domains.
\end{abstract}

\maketitle

It\^o's formula is one of the main tools in Stochastic
Analysis and, in particular, in the theory of
stochastic partial differential equations 
(SPDEs) of It\^o type. E. Pardoux (\cite{Pa75}) was the first to
consider the most general SPDEs with deterministic and stochastic
terms containing the unknown function and its derivatives
from an abstract point of view
of stochastic It\^o equations in the setting symbolically
described as $V\subset H\subset V^{*}$, where $H$
is a Hilbert space and $V,V^{*}$ are Banach spaces
(see \cite{Pa75} for references to previous results).
 One of the main steps in treating SPDEs consists of
establishing It\^o's formula for the square of the
 $H$-norm of solutions.

In   the deterministic case (without any stochastic terms)
 we deal with a function $v_{t}$, that is
in $V$ for almost all $t>0$ and its time derivative is 
in $V^{*}$, which is dual to $V$, for almost all $t>0$
(cf. Remark \ref{remark 8.13.1}).
 The goal is to show that
there is a modification $u_{t}$
 of $v_{t}$ which is an $H$-valued
{\em continuous\/} function and $\|u_{t}\|^{2}_{H}$
is an absolutely continuous function  admitting a natural
formula for its time derivative.
Even in this case  and even if $V,V^{*}$ are Hilbert spaces
the formula is not completely trivial.
For instance,
the proof of Theorem 3 on page 287 in \cite{Ev}
still has a tiny gap since the continuity of
$\|u_{t}\|^{2}_{H}$ at $t=0$ is not proved. 

In the stochastic case the proof given in
\cite{Pa75} and \cite{Pa79} 
is rather involved and consists of many steps.
In particular, it is based on the deterministic version
of It\^o's formula in Banach spaces with reference to
\cite{Li} where in Remark 1.2 on page 156 and 
Remark 7.9 on page 236 one indeed finds the statement
that the result is true but   neither
a proof or a reference to a proof is given.  
In \cite{KR} an  approach not using the deterministic
result  was suggested
for equations driven by continuous martingales.
In contrast with the deterministic case or, for that matter,
with \cite{Pa75} and \cite{Pa79}, where the
deterministic case was the starting point, the method of 
\cite{KR} is based on discretization of the time variable
and some arithmetical manipulations.
The method of \cite{KR} is generalized in \cite{GK} for the case
of arbitrary cadlag martingales, which required
 a dramatic 
increase in what concerns arithmetics.

Most likely there is no simple proof of It\^o's formula
in the case of general Banach spaces and this makes
it hard for a person interested in SPDEs to enter
the area. On the other hand, the case that $V,H$, and $V'$
are Hilbert spaces is the most common in applications
 and it turns out that
in this case there is a simple proof presented here
of
the continuity of the process $\|u_{t}\|^{2}_{H}$ and of
It\^o's formula for it.
The method
seems to be much simpler than the ones
 previously known. 

Basically, we lift $du_{t}$
  into $H$ by using resolvent operators
$R_{\lambda}$ constructed in Section \ref{section 7.17.1},
write It\^o's formula for Hilbert space-valued processes,
and then pass to the limit as $\lambda\to\infty$
in Section \ref{section 7.17.2}.
In the end of Section \ref{section 7.17.2} we give a version
of Theorem \ref{theorem 4.9.3} for SPDEs when $V$ is a
Sobolev-Hilbert space with integral
numbers of derivatives. In Section \ref{section 8.3.1}
we prove It\^o's formula not for the squared norm in $H$
but for more general functions as in \cite{Pa75}, \cite{Pa79}.
In our view the proofs given here
are more straightforward than previously known.
In the final Section \ref{section 8.6.1} we apply the result of
Section \ref{section 8.3.1} to prove the maximum principle
for the second-order SPDEs in arbitrary domains.
To the best of the author's knowledge,
in what concerns the arbitrariness of the domain
and the structure of the equation
 this is the most general result
known so far.  

It seems that the maximum principle 
for general second order SPDEs was first
proved in \cite{KR81} (see also \cite{KR82}
for the case of random coefficients) for SPDEs
in the whole space by the method of random characteristics
introduced there and also in   \cite{Ku} (a particular case
of the maximum principle appeared already in \cite{Pa75}).
Later the method of random characteristics was used
in many papers for various purposes, for instance,
to prove smoothness of solutions
(see, for instance, \cite{BG1}, \cite{BG2},
\cite{PT}, \cite{Tu} and the references therein). 
It was very tempting
to try to use this method for proving the maximum
principle for SPDEs in domains.
However, the implementation of the method
 turns out to become extremely
cumbersome and
inconvenient if the coefficients of the equation
are random processes. Also, it  requires
more regularity of solutions than actually needed.

In \cite{Kr07} a new method was suggested based on It\^o's
formula for $\|u^{+}_{t}\|^{2}_{H}$ which was derived
 by using mollifications in space variable. 
This required the underlying domain
to be in $C^{1}_{loc}$. In the present article the domain
is an arbitrary open set. On the other hand, the results
of \cite{Kr07} are much more general in other respects.
In particular, they were applied in \cite{Kr07} to prove the H\"older
continuity of solutions up to the boundary and the results
of the present article seem not to be applicable
for this purpose. Still they can be applied
in the same way as in \cite{KW} for
investigating filtering problems in the situation
of partially observable diffusion processes when the observations
are only available until the unobservable component
exits from a given domain. This would show that the result
of \cite{KW} 
 about filtering density before the exit time
occurs is true in case of arbitrary domains.

Finally, we mention that there are many situations
in which It\^o's formula is known for Banach space valued processes.
See, for instance, \cite{Br} and the references therein.
These formulas could be more general in some respects but they
do not cover our situation when the stochastic differential
leaves in $V^{*}$.

\mysection{Resolvent operator}
                                        \label{section 7.17.1}

It is highly unlikely that the results of this section
are new. However, the author's several attempts
to find them in the literature failed and were abandoned
for the reason that it takes less time 
  to prove the results
than to find them published somewhere. In addition,
the proof only takes two pages.

Let $V$ and $H$ be two Hilbert spaces with 
scalar products and norms $(\cdot,\cdot)_{V}$, $\|\cdot\|_{V}$
and $(\cdot,\cdot)_{H}$, $\|\cdot\|_{H}$, respectively. 
Assume that
$V\subset H$ and $V$ is dense in $H$ (in the metric of $H$)
and $\|u\|_{H}\leq \|u\|_{V}$ for any $u\in V$.

The  norm in $V$ is obviously equivalent to
$$
\big(\lambda\|u\|^{2}_{H}+\|u\|^{2}_{V}\big)^{1/2},
$$
where $\lambda\geq0$ is any fixed number.  
Then take an $f\in H$ and observe that the linear functional
$ (f,u)_{H}$ is bounded as a linear functional on $V$.
By Riesz's representation theorem there exists
a unique $v=:R_{\lambda}f\in V$
such that
$$
 (f,u)_{H}=\lambda (v,u)_{H}
+(v,u)_{V}
\quad\forall u\in V,
$$
that is
\begin{equation}
                                                   \label{4.13.2}
 (f,u)_{H}=\lambda (  R_{\lambda}f,u)_{H}
+( R_{\lambda}f,u)_{V},
\quad\forall u\in V,
\end{equation}
\begin{equation}
                                               \label{7.27.3}
( R_{\lambda}f,u)_{V}=((1-\lambda R_{\lambda})f,u )_{H},
\quad\forall u\in V.
\end{equation}

If $ f,g\in H$, then \eqref{4.13.2} with $u= R_{\lambda}g$
reads
$$
 (f, R_{\lambda}g)_{H}=
\lambda( R_{\lambda}f,  R_{\lambda}g)_{H}
+( R_{\lambda}f, R_{\lambda}g)_{V},
$$
where the right-hand side is symmetric in $f,g$. So is the
left-hand side implying that $ R_{\lambda}$
is a symmetric operator in $H$. After that
\eqref{4.13.2}  shows that $( R_{\lambda}f,u)_{V}=
(f, R_{\lambda}u)_{V}$ if $f,u\in V$,
so that $  R_{\lambda}$ is also a symmetric operator in 
$V$.

Then observe that
for $ u= R_{\lambda}f$ equation \eqref{4.13.2} implies that
\begin{equation}
                                                   \label{7.27.1}
\lambda\|  R_{\lambda}f\|^{2}_{H}+
 \| R_{\lambda}f\|^{2}_{V}
= (f, R_{\lambda}f)_{H}\leq  \|f\|_{H}\|  R_{\lambda}f
\|_{H},
\end{equation}
which yields the energy
estimates
\begin{equation}
                                                   \label{6.18.2}
 \lambda \| R_{\lambda}f\| _{H}\leq
\|f\|_{H},\quad\|R_{\lambda}f\|_{V}\leq\|f\|_{H}.
\end{equation}

\begin{theorem}
                                          \label{theorem 4.13.1}
(i) The norms of the operator $ \lambda R_{\lambda}$ as an operator from 
$H$ into $H$
 as well as an operator from $V$ into $V$
 are less than one;

(ii) If $f\in H$, $\lambda\geq0$, and $\lambda R_{\lambda}
f=f$, then $f=0$;

(iii) The set $ R_{\lambda}H$ is dense in 
$V$ in the metric of $V$;

(iv) For any $f\in H $ we have
\begin{equation}
                                                   \label{4.9.5}
\lim_{\lambda\to\infty}\|f-\lambda R_{\lambda}f\|_{H}=0;
\end{equation}

(v) For
$f\in V$ we have
\begin{equation}
                                                   \label{4.9.3}
\lim_{\lambda\to\infty}\|f-\lambda R_{\lambda}f\|_{V}=0.
\end{equation}
\end{theorem}
 
Proof.
(i) We get the first part of the assertion from
 \eqref{6.18.2}.
Next, if $f\in V$, then  
 for $u=f$
 we get from \eqref{7.27.3}
that
$$
( R_{\lambda}f,f)_{V}=((1-\lambda R_{\lambda})f,
(1-\lambda R_{\lambda})f)_{H}+((1-\lambda R_{\lambda})f,
\lambda R_{\lambda})f)_{H}
$$
$$
= \|f-\lambda R_{\lambda}f\|_{H}^{2}+
 (\lambda R_{\lambda}f,f)_{H}- \|\lambda R_{\lambda}f
\|^{2}_{H},
$$
where according to \eqref{7.27.1}
we have $ (\lambda R_{\lambda}f,f)_{H}- \|\lambda R_{\lambda}f
\|^{2}_{H}=\lambda \|R_{\lambda}f\|^{2}_{V}$, so that
\begin{equation}
                                         \label{7.27.2}
(  R_{\lambda}f,f)_{V}= \|f-\lambda R_{\lambda}f\|_{H}^{2}+
\lambda\| R_{\lambda}f\|^{2}_{V},
\end{equation}
$$
\lambda \|R_{\lambda}f\|^{2}_{V}\leq ( R_{\lambda}f,f)_{V},
\quad \|\lambda R_{\lambda}f\| _{V}\leq \| f\| _{V}.
$$
 
(ii) Under given conditions we have $f\in V$ and
equation \eqref{4.13.2} implies that $(R_{\lambda} f,u)_{V}=0$
for all $u\in V$. Hence $R_{\lambda} f=0$ and $f=0$ indeed.

(iii) Assume the contrary. Then there exists $u\in V$,
$u\ne0$, such that
$(  R_{\lambda}f,u)_{V}=0$ for all $f\in H$. Then \eqref{4.13.2}
shows that $(f,u)_{H}=(\lambda R_{\lambda}f,u)_{H}
=(f,\lambda R_{\lambda}u)_{H}$ for all $f\in H$. It follows that
$\lambda R_{\lambda}u=u$ and $u=0$ by (ii), which is the desired 
contradiction.

(iv) If $f\in V$, this assertion follows from \eqref{7.27.2} 
after we   let $\lambda\to\infty$,
and use (i). In the general case it suffices to use
the denseness of $V$ in $H$ and assertion (i).

(v) Owing to (iii) and (i) while proving \eqref{4.9.3}
we may concentrate  on $f=R_{1}g$, where $g\in H$.
Next, we observe that taking $u=(1-\lambda R_{\lambda})f$
in \eqref{7.27.3} leads to
$$
( R_{\lambda}f,(1-\lambda R_{\lambda})f)_{V}=((1-\lambda R_{\lambda})f
, (1-\lambda R_{\lambda})f)_{H},
$$
which implies that
$$
\|(1-\lambda R_{\lambda})f\|_{V}^{2}=
(f,(1-\lambda R_{\lambda})f)_{V}-(\lambda R_{\lambda}f,(1-\lambda R_{\lambda})f)_{V}
$$
$$
=(f,(1-\lambda R_{\lambda})f)_{V}-\lambda\|(1-\lambda R_{\lambda})f
\|^{2}_{H}\leq (f,(1-\lambda R_{\lambda})f)_{V}.
$$
Here, in light of \eqref{7.27.3}, the last expression is
$$
(R_{1}g,(1-\lambda R_{\lambda})f)_{V}=
((1-R_{1})g,(1-\lambda R_{\lambda})f )_{H},
$$
which tends to zero as $\lambda\to\infty$ by (iv).
The theorem is proved. 

As a justification of the notation $R_{\lambda}$ and its name
as a resolvent operator consider the following
situation.

Let $G$ be an open set in $\bR^{d}=\{x=(x^{1},...,x^{d}):
x^{i}\in\bR\}$ and let $m\geq1$ be an integer.
Define
$$
\HO^{m}_{2}=\HO^{m}_{2}(G)
$$
as the closure of $C^{\infty}_{0}=
C^{\infty}_{0}(G)$ with respect to the norm
$$
\|u\|_{H^{m}_{2}}:=
\big(\sum_{|\alpha|\leq m}\int_{G} 
|D^{\alpha}u|^{2}\,dx\big)^{1/2},
$$
where as usual for any multi-index $\alpha=(\alpha_{1},
...,\alpha_{d})$
$$
D^{\alpha}=D^{\alpha_{1}}_{1}\cdot...\cdot D^{\alpha_{d}}_{d},
\quad D_{i}=\partial/\partial x^{i},\quad
|\alpha|=\alpha_{1}+...+\alpha_{d}.
$$
The space $\HO^{m}_{2}$ is a Hilbert space with scalar product
$$
(u,\phi)_{H^{m}_{2}}=\sum_{|\alpha|\leq m}
\int_{G}(D^{\alpha}u)D^{\alpha}\phi\,dx.
$$

If above we take $H=L_{2}=L_{2}(G)$ and $V=\HO^{m}_{2}$
then our hypotheses about $V$ and $H$ are satisfied and
\eqref{4.13.2} becomes
$$
\int_{G}fu\,dx=\lambda\int_{G}uR_{\lambda}f\,dx
+\sum_{|\alpha|\leq m}
\int_{G}(D^{\alpha}R_{\lambda}f)D^{\alpha}u\,dx,
$$
which in the sense of generalized functions shows that
$R_{\lambda}f$ is a solution of the equation
$$
f=\lambda v-Lv,
$$
where
$$
L=-\sum_{|\alpha|\leq m}(-1)^{|\alpha|}D^{2\alpha}.
$$
Hence, $R_{\lambda}$ is indeed a resolvent operator for $L$.
By the way, this example shows that, generally, $R_{\lambda}H\ne V$.

\mysection{It\^o's formula for the squared norm}
                                       \label{section 7.17.2}

Let $(\Omega,\cF,P)$ be a complete probability space
and let $\{\cF_{t},t\geq\}$ be an increasing filtration
of $\sigma$-fields $\cF_{t}\subset \cF$, which are complete
with respect to $\cF,P$.  

In order to avoid unimportant complications
we assume that $(V,(\cdot,\cdot)_{V})$ is 
a separable Hilbert space, which is the case
in many applications. Then $(H,(\cdot,\cdot)_{H})$
is also separable. It is convenient that
under this assumption
there is no difference between weak and strong measurability.

Assume that we are given $V $-valued
processes $v_{t},v^{*}_{t}$, $t>0$,
 which are predictable and satisfy
\begin{equation}
                                                   \label{8.9.5}
E\int_{0}^{T}\|v_{t},v^{*}_{t}\|^{2}_{V }
\,dt<\infty
\end{equation}
  for any $T\in(0,\infty)$. Also let  $m_{t}$, $t\geq0$,
be an $H $-valued
continuous martingale starting at the origin with
\begin{equation}
                                                     \label{8.9.1}
d\langle m\rangle_{t}\leq dt.
\end{equation}

The theory of integrating predictable 
Hilbert-space valued processes with respect to
 continuous same space-valued martingales is quite
parallel to that in case the Hilbert space is just $\bR^{d}$.
This theory implies that under the above conditions 
the stochastic integral
\begin{equation}
                                               \label{6.18.6}
h_{t}:=\int_{0}^{t}(v_{s},dm_{s})_{H }
\end{equation}
is well defined and  is a
continuous real-valued martingale with
$$
\langle h\rangle_{t}\leq \int_{0}^{t}\|v_{s}\|^{2}_{H }\,ds.
$$
Suppose that $v_{0}$ is
an 
$H $-valued $\cF_{0}$-measurable random vector.
Finally, assume that for any $\phi\in V$
we have
\begin{equation}
                                                 \label{4.9.7}
(\phi,v_{t})_{H }=(\phi,v_{0})_{H }+\int_{0}^{t}
(\phi,v^{*}_{s})_{V }\,ds+(\phi,m_{t})_{H }
\end{equation}
for almost all $(\omega,t)$.

\begin{theorem}
                                             \label{theorem 4.9.3}
Under the above assumptions there exists a continuous
$H $-valued $\cF_{t}$-adapted process $u_{t}$
and a set $\Omega'\subset\Omega$ of full probability
such that 

(i) $u_{t}=v_{t}$ for almost all $(\omega,t)$, so that
$$
E\int_{0}^{T}\|u_{t} \|^{2}_{V }
\,dt<\infty
$$
  for any $T\in(0,\infty)$,

(ii) for
all $\omega\in\Omega'$, all $\phi\in V$,
and all $t\geq0$ we have
\begin{equation} 
                                                 \label{4.9.8}
(\phi,u_{t})_{H }=(\phi,v_{0})_{H }+\int_{0}^{t}
(\phi,v^{*}_{s})_{V }\,ds+(\phi,m_{t})_{H },
\end{equation}

(iii) for
all $\omega\in\Omega'$ 
and all $t\geq0$ we have 
\begin{equation}
                                                 \label{4.9.9}
\|u_{t}\|_{H }=\|v_{0}\|_{H }+2\int_{0}^{t}
(u_{s},v^{*}_{s})_{V }\,ds
+\langle m\rangle_{t}
+2 \int_{0}^{t}(v_{s},dm_{s})_{H }.
\end{equation}

\end{theorem}

Proof. Inspired by \eqref{7.27.3}
for $n=1,2,...$ define $S_{n}=nR_{n}$ and
\begin{equation}
                                                    \label{8.6.1}
u^{n}_{t}=S_{n}v_{0}+\int_{0}^{t}n (1-S_{n})  
v^{*}_{s}\,ds+ S_{n}m_{t}.
\end{equation}
Here the integral  makes sense as the integral  of an $H $-valued
function.
Furthermore, $u^{n}_{t}$ is obviously continuous
as an $H $-valued function.

Also observe that \eqref{4.9.7}  with $\phi=S_{n}\psi$,
$\psi\in H $, and \eqref{7.27.3} yield that for almost all
$(\omega,t)$
$$
(\psi,S_{n}v_{t})_{H }=(\psi,S_{n}v_{0})_{H }
+\int_{0}^{t}(\psi, n(1-S_{n}) 
v^{*}_{s})_{H }\,ds
$$
\begin{equation}
                                            \label{8.6.2}
+(\psi, S_{n}m_{t})_{H }=
(\psi,u^{n}_{t})_{H }.
\end{equation}
This and the separability of $H $ shows that 
$$
u^{n}_{t}=S_{n}v_{t}
$$ 
for almost all $(\omega,t)$.

 Next, from Doob's inequality
it follows that for any $T\in[0,\infty)$
$$
 E\sup_{t\leq T}\|u^{n}_{t} \|^{2}_{H }<\infty.
$$

By It\^o's formula for integrals of Hilbert-space valued processes
we have (a.s.)
$$
\|u^{n}_{t}\|_{H }^{2}=\|S_{n}v_{0}\|_{H }^{2}
+2\int_{0}^{t}(u^{n}_{s}, n(1-S_{n}) 
v^{*}_{s})_{H } \,ds
$$
\begin{equation}
                                                    \label{4.11.1}
+\langle S_{n}m\rangle_{t}
+2\int_{0}^{t}(S_{n}u^{n}_{s}, dm_{s})_{H },
\end{equation}
$$
\|u^{n}_{t}-u^{k}_{t}\|^{2}_{H }=\|(S_{n}-S_{k})v_{0}\|_{H }^{2}
+\langle (S_{n}-S_{k})m\rangle_{t}
+2\int_{0}^{t}(S_{n}u^{n}_{s}-S_{k}u^{k}_{s}, dm_{s})_{H }
$$
\begin{equation}
                                                    \label{4.11.2}
+2\int_{0}^{t}(u^{n}_{s}-u^{k}_{s},
 [n(1-S_{n}) -k(1-S_{k})  ]
v^{*}_{s})_{H } \,ds
\end{equation}
for all $t\geq0$.

Observe that there is $\Omega'$ with $P(\Omega')=1$ such that
\begin{equation}
                                                 \label{6.19.5}
u^{n}_{t}=S_{n}v_{t},\quad n=1,2,...,\quad
 v_{t},v^{*}_{t}\in V 
\end{equation}
 for almost all $t$ on $\Omega'$. It follows that
in the integrands in \eqref{4.11.1} and \eqref{4.11.2}
 we can replace
$u^{n}_{s}$ with $S_{n}v_{s}$ if $\omega\in\Omega'$
and use \eqref{7.27.3}. Then
for $\omega\in\Omega'$ and $t$ such that \eqref{6.19.5} holds
we have
$$
(u^{n}_{s}-u^{k}_{s},
 [n(1-S_{n}) -k(1-S_{k}) ]
v^{*}_{s})_{H }=
$$

$$
(S_{n}v_{s}-S_{k}v_{s},
 n(1-S_{n}) v^{*}_{s})_{H }
-(S_{n}v_{s}-S_{k}v_{s},k(1-S_{k}) 
v^{*}_{s})_{H }
$$

$$
=(S_{n}v_{s}-S_{k}v_{s},S_{n}v^{*}_{s})_{V }-
(S_{n}v_{s}-S_{k}v_{s},S_{k}v^{*}_{s})_{V }
$$

$$
=(S_{n}v_{s}-S_{k}v_{s},S_{n}v^{*}_{s}-
S_{k}v^{*}_{s})_{V }.
$$

Hence, for $\omega\in\Omega'$ and all $n,k\geq1$ and
 $t\geq0$ we get that
$$
\|u^{n}_{t}\|_{H }^{2}=\|S_{n}v_{0}\|_{H }^{2}
+2\int_{0}^{t}(S_{n}v_{s},
v^{*}_{s})_{V } \,ds
$$
$$
+\langle S_{n}m\rangle_{t}
+2\int_{0}^{t}(S_{n}u^{n}_{s}, dm_{s})_{H },
$$
$$
\|u^{n}_{t}-u^{k}_{t}\|^{2}_{H }=\|(S_{n}-S_{k})v_{0}\|_{H }^{2}
+\langle (S_{n}-S_{k})m\rangle_{t}
+2\int_{0}^{t}(S^{2}_{n}v_{s}-S^{2}_{k}v_{s}, dm_{s})_{H }
$$
$$
+2\int_{0}^{t}(S_{n}v_{s}-S_{k}v_{s},S_{n}v^{*}_{s}-
S_{k}v^{*}_{s})_{V } \,ds.
$$

Furthermore, by Doob's inequality for any $T\in[0,\infty)$
\begin{equation}
                                                    \label{4.11.4}
E\sup_{t\leq T}\|u^{n}_{t}-u^{k}_{t}\|^{2}_{H }
\leq I_{nk}^{1}+2I^{2}_{nk}+I^{3}_{nk}+4(I^{4}_{nk})^{1/2},
\end{equation}
where
$$
I_{nk}^{1}=E\|(S_{n}-S_{k})v_{0}\|_{H }^{2}
$$
$$
I_{nk}^{2}= E\int_{0}^{T}|((S_{n}-S_{k})v_{s},
(S_{n}-S_{k})
v^{*}_{s})_{V} |\,ds,
$$
$$
I_{nk}^{3}=E
\langle (S_{n}-S_{k})m\rangle_{T}=E\|(S_{n}-S_{k})m_{T}\|^{2}
_{H },
$$
$$
I_{nk}^{4}=E \int_{0}^{T} \|S^{2}_{n}v_{s}
-S^{2}_{k}v_{s}\|_{H }^{2}\,ds .
$$
By using the dominated convergence theorem, 
Theorem \ref{theorem 4.13.1}, and the inequality
$$
|((S_{n}-S_{k})v_{s}, (S_{n}-S_{k})
v^{*}_{s})_{V } |
$$
$$
\leq\|(S_{n}-S_{k})v_{s}\|^{2}_{V }
+\|(S_{n}-S_{k})v^{*}_{s}\|^{2}_{V },
$$
we easily conclude that $I_{nk}^{1}+2I^{2}_{nk}+I^{3}_{nk}\to0$
as $n,k\to\infty$. Furthermore, $S_{n}^{2}-S_{k}^{2}
=(S_{n}+S_{k})(S_{n}-S_{k})$ so that
$$
I^{4}_{nk}\leq 4 E\int_{0}^{T}\|(S_{n}-S_{k})
v_{s}\|_{H }^{2} \,ds, 
$$
which by the dominated convergence theorem
implies that $I^{4}_{nk}\to0$ 
as $n,k\to\infty$ as well.

 We now conclude from \eqref{4.11.4} that its left-hand side
tends to zero. Furthermore,
$$
E\int_{0}^{T}\|u^{n}_{t}-v_{t}\|^{2}_{V }\,dt
=E\int_{0}^{T}\|(S_{n} -1)v_{t}\|^{2}_{V }\,dt\to0.
$$
Hence $u^{n}_{t}$ converges to $v_{t}$ in $H 
(\Omega\times(0,T),V)$
and converges uniformly on $[0,T]$ as $H $-valued functions
in probability. The latter limit we denote by $u_{t}$
and show that this function is the one we want.
Of course, $u_{t}$ is a continuous $H $-valued functions,
 it is $\cF_{t}$-adapted, and $u_{t}=v_{t}$ for almost all
$(\omega,t)$.

One easily obtains that for each $t$ 
equation \eqref{4.9.9} holds with probability
one by passing to the limit in \eqref{4.11.1}.
Since both parts of \eqref{4.9.9} are continuous in $t$,
it holds on a set of full probability for all $t$.

Obviously  \eqref{4.9.7} will hold for almost all $(\omega,t)$
if we replace $v_{t}$ with $u_{t}$, that is,  \eqref{4.9.8}
holds for any $\phi\in V$ for almost all
$(\omega,t)$. The continuity of both parts of 
\eqref{4.9.8} with respect to $t$ and $\phi\in V$
and the separability of $V$ then imply that there is
a set $\Omega'$ of full probability such that assertion (iii)
holds.

The theorem is proved.
\begin{remark}
The reader understands, of course, that condition
\eqref{8.9.5} can be  replaced with the same condition but
without expectation sign. This generalization
is easily achieved by using appropriate stopping times.
\end{remark}

Next in the setting described in the end of Section \ref{section 7.17.1}
suppose that we are given an $\HO^{m}_{2}$-valued
process $v_{t}$ and $L_{2}$-valued processes
$f^{\alpha}_{t} $, $|\alpha|\leq m$. We assume that all these processes
are predictable and such that
$$
E\int_{0}^{T}\big[\|v_{t}\|^{2}_{H^{m}_{2}}+
\sum_{|\alpha|\leq m}\|f^{\alpha}_{t}\|^{2}_{L_{2}}\big]\,dt<\infty
$$
for any $T\in(0,\infty)$. We also assume that we are
given a continuous $L_{2}$-valued martingale $m_{t}$
satisfying \eqref{8.9.1}
and   $v_{0}$ is
an 
$L_{2}$-valued $\cF_{0}$-measurable random function.
 Finally, suppose that
for any $\phi\in C^{\infty}_{0}$ we have
\begin{equation}
                                                   \label{8.9.2}
\int_{G}\phi v_{t}\,dx=\int_{G}\phi v_{0}\,dx
+\int_{0}^{t}\int_{G} \sum_{|\alpha|\leq m}f^{\alpha}_{s}
D^{\alpha}\phi
\,dx\,ds
+\int_{G}\phi m_{t}\,dx
\end{equation}
for almost all $(\omega,t)$.
\begin{remark}
                                           \label{remark 8.13.1}
Formally \eqref{8.9.2} can be expressed as
$$
dv_{t}=\sum_{|\alpha|\leq m}(-1)^{|\alpha|}D^{\alpha}
f^{\alpha}_{t}\,dt + dm_{t},
$$
which shows that $dv_{t}$ lives in $H^{-m}_{2}
:=(\HO^{m}_{2})^{*}$.
\end{remark}
  
\begin{theorem}
                                             \label{theorem 8.9.1}
Under the above assumptions there exists a continuous
$L_{2}$-valued $\cF_{t}$-adapted process $u_{t}$
and a set $\Omega'\subset\Omega$ of full probability
such that 

(i) $u_{t}=v_{t}$ for almost all $(\omega,t)$, so that
$$
E\int_{0}^{T}\|u_{t} \|^{2}_{H^{m}_{2} }
\,dt<\infty
$$
  for any $T\in(0,\infty)$,

(ii) for
all $\omega\in\Omega'$, all $\phi\in \HO^{m}_{2}$,
and all $t\geq0$ we have
\begin{equation} 
                                                 \label{4.9.08}
(\phi,u_{t})_{L_{2}}=(\phi,v_{0})_{L_{2}}+\int_{0}^{t}
\sum_{|\alpha|\leq m}
(D^{\alpha}\phi,f^{\alpha}_{s})_{L_{2}}\,ds+(\phi,m_{t})_{L_{2}},
\end{equation}

(iii) for
all $\omega\in\Omega'$ 
and all $t\geq0$ we have 
$$
\|u_{t}\|_{L_{2}}=\|v_{0}\|_{L_{2}}+2\int_{0}^{t}
\sum_{|\alpha|\leq m}
(D^{\alpha}u_{s},f^{\alpha}_{s})_{L_{2}}\,ds
+\langle m\rangle_{t}
+2 \int_{0}^{t}(v_{s},dm_{s})_{L_{2}}.
$$

\end{theorem}

Proof. To derive this result from Theorem
\ref{theorem 4.9.3}, we first observe that
in light of the denseness of $C^{\infty}_{0}$
in $\HO^{m}_{2}$ equation \eqref{8.9.2} also holds for any
$\phi\in\HO^{m}_{2}$. Then notice that for each $(\omega,s)$
$$
F_{s}(\phi):=\int_{G}\sum_{|\alpha|\leq m}f^{\alpha}_{s}
D^{\alpha}\phi\,dx
$$
is a bounded linear functional on  $\HO^{m}_{2}$
with
$$
|F_{s}(\phi)|\leq\|\phi\|_{H^{m}_{2}}\big(
\sum_{|\alpha|\leq m}\|f^{\alpha}_{s}\|_{L_{2}}^{2}\big)^{1/2}.
$$
It follows by Riesz's representation theorem that there
exists a unique $v^{*}_{s}\in\HO^{m}_{2}$ such that 
$$
F_{s}(\phi)=(\phi,v^{*}_{s})_{H^{m}_{2}},\quad
\|v^{*}_{s}\|^{2}\leq
\sum_{|\alpha|\leq m}\|f^{\alpha}_{s}\|_{L_{2}}^{2}.
$$
These relations imply that $v^{*}_{s}$ is weakly 
predictable and, since $\HO^{m}_{2}$ is separable, it
is (just) predictable. Also we have that condition
\eqref{8.9.5} is satisfied. Hence one can rewrite
\eqref{8.9.2}
in form \eqref{4.9.7} and then all assertion
of the present theorem follow directly from
Theorem \ref{theorem 4.9.3}. The theorem is proved.

\mysection{A more general It\^o's formula}
                                         \label{section 8.3.1}
We suppose that all assumptions stated in 
Section \ref{section 7.17.2} are satisfied.

Let $\phi(h)$ be a real-valued function on $H$.
Assume that 

(i) for any $h,\xi\in H$   the functions 
$\phi(h+t\xi)$ is   twice continuously differentiable
as a function of $t$ and the functions
$$
\phi_{(\xi)}(h):=\frac{\partial}{\partial t}\phi(h+t\xi)_{t=0},
\quad
\phi_{(\xi)(\xi)}(h):=\frac{\partial^{2}}{
(\partial t)^{2}}\phi(h+t\xi)\big|_{t=0}
$$
are continuous as functions of $(h,\xi)\in H\times H$;

(ii) For any $R\in(0,\infty)$
there exists a $K(R)$ such that for all 
$h,\xi\in H$  satisfying $\|h\|_{H}\leq R$
we have
$$
|\phi_{(\xi)}(h)|\leq K(R)\|\xi\|_{H},\quad
|\phi_{(\xi)(\xi)}(h)|\leq K(R)\|\xi\|^{2}_{H}.
$$
In this situation for any $h\in H$ the function
$\phi_{(\xi)}(h)$ as a function of $\xi\in H$
is a continuous linear functional and by Riesz's
representation theorem there exists 
an element $\phi_{(\cdot)}(h)\in H$ such that
$$
\phi_{(\xi)}(h)=(\phi_{(\cdot)}(h),\xi)_{H},\quad
\|\phi_{(\cdot)}(h)\|_{H}\leq K(\|h\|_{H}).
$$
Next, we assume that, 

(iii) If $h\in V$, then $\phi_{(\cdot)}(h)\in V$
and
$$
\|\phi_{(\cdot)}(h)\|_{V}\leq K(1+\|h\|_{V}),
$$
where $K$ is a fixed constant;

(iv) For any $v^{*}\in V$ the function $(\phi_{(\cdot)}(v),v^{*})_{V}$
is a continuous function on $V$ (in the metric of $V$).

Let $w^{1}_{t},w^{2}_{t},...$ be a finite or infinite sequence
of independent Wiener processes on $(\Omega,\cF,P)$,
which are Wiener processes with respect to $\{\cF_{t}\}$.
We assume that we are given a sequence of predictable
$H$-valued processes $\sigma^{i}_{t}$ such that
for any $T\in(0,\infty)$  
$$
\sum_{k}E\int_{0}^{T}\|\sigma^{k}_{t}\|_{H}^{2}\,dt<\infty.
$$
Under this assumption it is well known that the series
$$
\sum_{k}\int_{0}^{t}\sigma^{k}_{s}\,dw^{k}_{s}
$$
converges in $H$ uniformly on finite time intervals
in probability and we assume that the series converges
to $m_{t}$. From the continuity of the scalar product
in $H$ it follows also that for any $h\in H$ we have
(a.s.) for all $t$
$$
(h,m_{t})_{H}=\sum_{k}\int_{0}^{t}(h,\sigma^{k}_{s})_{H}\,dw^{k}_{s},
$$
where the series converges  uniformly on finite time intervals
in probability. Then equation \eqref{4.9.08} is equivalent
to saying
that the function $u_{t}$
 satisfies
\begin{equation}
                                                 \label{8.7.1}
(\phi,u_{t})_{H }=(\phi,v_{0})_{H }+\int_{0}^{t}
(\phi,v^{*}_{s})_{V }\,ds+
\sum_{k}\int_{0}^{t}(\phi,\sigma^{k}_{s})_{H}\,dw^{k}_{s}
\end{equation}
for each $\phi\in V$ (a.s.) for all $t$.

The following result can be found in \cite{Pa75}
in a more general situation. Our innovation is
a different and   shorter proof.
\begin{theorem}
                                          \label{theorem 8.3.1}
Under the above assumptions (a.s.) for all $t$
$$
\phi(u_{t})=\phi(u_{0})+
\sum_{k}\int_{0}^{t} \phi_{(\sigma^{k}_{s})}(u_{s}) \,dw^{k}_{s}
$$
\begin{equation}
                                             \label{8.3.2}
+\int_{0}^{t} \big[(\phi_{(\cdot)}(u_{s}), v^{*}_{s})_{V} 
+(1/2)\sum_{k} \phi_{(\sigma^{k}_{s})
(\sigma^{k}_{s})}(u_{s})\big]\,ds,
\end{equation}
where $u_{t}$ is taken from Theorem \ref{theorem 4.9.3}
and the series of stochastic integrals
 converges uniformly on finite time intervals
in probability.
\end{theorem}

Proof. The last assertion of the theorem follows from the
fact that the series of quadratic variations
of the stochastic integrals in \eqref{8.3.2} converges:
$$
\sum_{k}\int_{0}^{t}|\phi_{(\sigma^{k}_{s})}(u_{s})|^{2}
\,ds\leq\sum_{k}\int_{0}^{t}\|\phi_{(\cdot)}(u_{s})\|_{H}^{2}
\|\sigma^{k}_{s}\|_{H}^{2}\,ds
$$
$$
\leq K^{2}\big(\sup_{s\leq t}\|u_{s}\|_{H}\big)
\sum_{k}\int_{0}^{t} 
\|\sigma^{k}_{s}\|_{H}^{2}\,ds<\infty.
$$
It is also worth noting that other terms in \eqref{8.3.2}
make sense as well. Indeed,
$$
\int_{0}^{t} \sum_{k} |\phi_{(\sigma^{k}_{s})
(\sigma^{k}_{s})}(u_{s})|\,ds\leq
K^{2}\big(\sup_{s\leq t}\|u_{s}\|_{H}\big)
\sum_{k}\int_{0}^{T}\|\sigma^{k}_{t}\|_{H}^{2}\,dt<\infty,
$$
$$
\int_{0}^{t} |(\phi_{(\cdot)}(u_{s}), v^{*}_{s})_{V}|\,ds
\leq K\int_{0}^{t}(1+\|u_{s}\|_{V})\|v^{*}\|_{V}\,ds<\infty.
$$
This argument shows that the right-hand side of
\eqref{8.3.2} is a  continuous process (a.s.). So is
its left-hand side and, to prove that
\eqref{8.3.2} holds (a.s.) for all $t$,
it suffices to prove that \eqref{8.3.2} holds  for  
each $t$ (a.s.).

The rest of the proof we split into a few steps.

{\em Step 1\/}. Consider the case that $V=H$
and the number of the Wiener processes
is finite, say, equal to $p$. Take an orthonormal
 basis $\{e_{i}\}$ in $H$, denote by $\Pi_{n}$
the orthogonal projection operator on ${\rm Span}\,\{e_{1},...,e_{n}
\}$, and set
$$
u^{n}_{t}:=\Pi_{n}u_{t}=\Pi_{n}u_{0}+\int_{0}^{t}
\Pi_{n}v^{*}_{s}\,ds+\sum_{k\leq p}\int_{0}^{t}
\Pi_{n}\sigma^{k}_{s}\,dw^{k}_{s},\quad 
\phi_{n}(h)=\phi(\Pi_{n}h).
$$
The function $\phi_{n}$, as a 
continuous function on a finite-dimensional
Euclidean space,   has two continuous directional derivatives
in any direction. Therefore, it is twice continuously
differentiable and by the classical It\^o's formula
$$
\phi(u^{n}_{t})=\phi_{n}(u^{n}_{t})
=\phi(\Pi_{n}u_{0})+\sum_{k\leq p}\int_{0}^{t}
\phi_{(\Pi_{n}\sigma^{k}_{s})}(\Pi_{n}u _{s})\,dw^{k}_{s}
$$
\begin{equation}
                                                  \label{7.5.1}
+\int_{0}^{t}\big[\phi_{(\Pi_{n}v^{*}_{s})}
(\Pi_{n}u _{s})+(1/2)\sum_{k\leq p}\phi_{(\Pi_{n}\sigma^{k}_{s})
(\Pi_{n}\sigma^{k}_{s})}
(\Pi_{n}u _{s})\big]\,ds.
\end{equation}
Here
$$
|\phi_{(\Pi_{n}\sigma^{k}_{s})
(\Pi_{n}\sigma^{k}_{s})}
(\Pi_{n}u _{s})|\leq \|\sigma^{k}_{s}\|_{H}^{2}
K\big(\max_{s\leq t}\|u_{s}\|_{H}\big)
$$
and on an event of full probability on which $u_{s}$
is an $H$-valued continuous function
$$
\phi_{(\Pi_{n}\sigma^{k}_{s})
(\Pi_{n}\sigma^{k}_{s})}
(\Pi_{n}u _{s})\to \phi_{( \sigma^{k}_{s})
( \sigma^{k}_{s})}
( u _{s})
$$
for all $s$ and $k$. It follows by the dominated
convergence theorem that
$$
\int_{0}^{t}\phi_{(\Pi_{n}\sigma^{k}_{s})
(\Pi_{n}\sigma^{k}_{s})}
(\Pi_{n}u _{s}) \,ds\to \int_{0}^{t}\phi_{( \sigma^{k}_{s})
( \sigma^{k}_{s})}
(  u _{s}) \,ds
$$
for any $t$ and $k$ (a.s.). Similarly, for any $t$ (a.s.)
$$
\int_{0}^{t} \phi_{(\Pi_{n}v^{*}_{s})}
(\Pi_{n}u _{s})\,ds\to \int_{0}^{t} \phi_{( v^{*}_{s})}
( u _{s})\,ds.
$$
Finally, by the same reasons as above
$$
\int_{0}^{t}|\phi_{(\Pi_{n}\sigma^{k}_{s})}(\Pi_{n}u _{s})
-\phi_{( \sigma^{k}_{s})}( u _{s})|^{2}\,ds\to 0
$$
for any $t$ and $k$ (a.s.).

This allows us to pass to the limit in \eqref{7.5.1}
and conclude that
$$
\phi(u _{t}) 
=\phi( u_{0})+\sum_{k\leq p}\int_{0}^{t}
\phi_{( \sigma^{k}_{s})}( u _{s})\,dw^{k}_{s}
$$
\begin{equation}
                                                  \label{7.5.2}
+\int_{0}^{t}\big[\phi_{( v^{*}_{s})}
( u _{s})+(1/2)\sum_{k\leq p}\phi_{( \sigma^{k}_{s})
( \sigma^{k}_{s})}
( u _{s})\big]\,ds 
\end{equation}
for any $t$ (a.s.).

{\em Step 2\/}. Again let $V=H$ but suppose that the number
of the Wiener processes is infinite. Then introduce
$$
u^{n}_{t}=u_{0}+\int_{0}^{t}v^{*}_{s}\,ds+\sum_{k\leq n}
\int_{0}^{t}\sigma^{k}_{s}\,dw_{s}
$$
and observe that, as we pointed out before the theorem,
for any (finite) $t,\varepsilon>0$,
\begin{equation}
                                                  \label{7.5.4}
P(\sup_{s\leq t}\|u_{s}-u^{n}_{s}\|_{H}>\varepsilon)\to0
\end{equation}
as $n\to\infty$.
By the result of Step 1
$$
\phi(u^{n}_{t}) 
=\phi( u_{0})+\sum_{k\leq n}\int_{0}^{t}
\phi_{( \sigma^{k}_{s})}( u^{n} _{s})\,dw^{k}_{s}
$$
\begin{equation}
                                                  \label{7.5.3}
+\int_{0}^{t}\big[\phi_{( v^{*}_{s})}
( u^{n}_{s})+(1/2)\sum_{k\leq n}\phi_{( \sigma^{k}_{s})
( \sigma^{k}_{s})}
( u^{n}_{s})\big]\,ds 
\end{equation}
for any $t$ (a.s.).

Next, owing to \eqref{7.5.4} there is a subsequence $n(j)\to
\infty$ as $j\to\infty$ such that for any $t\in(0,\infty)$
(a.s.)
$$
\sup_{s\leq t}\|u_{s}-u^{n(j)}_{s}\|_{H}\to0.
$$
Then, of course,
$$
\int_{0}^{t} \phi_{( v^{*}_{s})}( u^{n(j)}_{s})\,ds
\to\int_{0}^{t} \phi_{( v^{*}_{s})}( u _{s})\,ds\quad
\text{(a.s.)}.
$$

Furthermore, (a.s.)
$$
\phi_{(\sigma^{k}_{s})}(u^{n(j)}_{s})\to\phi_{(\sigma^{k}_{s})}(u _{s})
$$
 because of the continuity
of $\phi_{(\xi)}(h)$ on $H\times H$. In addition,
$$
\sum_{k} |\phi_{\sigma^{k}_{s}}(u^{n(j)}_{s})|^{2}
 \leq K^{2}\big(\sup_{s\leq t,r\geq1}\|u^{n(r)}_{s}\|_{H}\big)
\sum_{k}\|\sigma^{k}_{s}\|_{H}^{2}\,ds
$$
and the right-hand side has a finite integral over $[0,t]$
(a.s.).
It follows by the dominated convergence theorem that
the quadratic variation at time $t$ of the difference
$$
\sum_{k\leq n}\int_{0}^{t}
\phi_{( \sigma^{k}_{s})}( u^{n(j)} _{s})\,dw^{k}_{s}-
\sum_{k}\int_{0}^{t}
\phi_{( \sigma^{k}_{s})}( u _{s})\,dw^{k}_{s}
$$
tends to zero (a.s.) as $k\to\infty$ and the difference itself
goes to zero in probability.

For similar reasons
$$
\sum_{k}\int_{0}^{t}|\phi_{( \sigma^{k}_{s})
( \sigma^{k}_{s})}
( u^{n(j)}_{s})-\phi_{( \sigma^{k}_{s})
( \sigma^{k}_{s})}
( u _{s}) | \,ds\to0 
$$
(a.s.) and we conclude from \eqref{7.5.3} that \eqref{8.3.2}
holds (a.s.).

{\em Step 3\/}. Now we consider the general case. As
in the proof of Theorem \ref{theorem 4.9.3} we introduce
$u^{n}_{t}$ by \eqref{8.6.1} 
and observe that the computation
\eqref{8.6.2} shows that (a.s.) $u^{n}_{t}=S_{n}u_{t}$
for all $t$. According to Step 2 for any $t$  (a.s.) 
$$
\phi(S_{n}u_{t})=\phi(S_{n}u_{0})
+\sum_{k}\int_{0}^{t} \phi_{(S_{n}\sigma_{s}^{i})}(S_{n}u_{s})
\,dw^{k}_{s}
$$
\begin{equation}
                                             \label{8.6.5}
+\int_{0}^{t}\big[(\phi_{(\cdot)}(S_{n}u_{s}),n(1-S_{n})v^{*}_{s})_{H}
+(1/2)\sum_{k}\phi_{(S_{n}\sigma^{k}_{s})(S_{n}\sigma^{k}_{s})}
(S_{n}u_{s})\big]
\,ds.
\end{equation}
Here (a.s.) for all $s$
$$
\phi_{(S_{n}\sigma^{k}_{s})}(S_{n}u _{s})\to
\phi_{(\sigma^{k}_{s})}(u _{s})
$$
as $n\to\infty$ because of the continuity
of $\phi_{(\xi)}(h)$ on $H\times H$. Furthermore,
$$
|\phi_{(S_{n}\sigma^{k}_{s})}(S_{n}u _{s})|\leq
K(\|u_{s}\|_{H})
\|S_{n}\sigma^{k}_{s}\|_{H}\leq
K(\|u_{s}\|_{H})
\| \sigma^{k}_{s}\|_{H}.
$$
As before this implies that the series of stochastic integrals in 
\eqref{8.6.5} converges to that in \eqref{8.3.2}
in probability as $n\to\infty$.

Next, owing to \eqref{7.27.3} and the fact that $S_{n}u_{s}\in V$
and $\phi_{(\cdot)}(S_{n}u_{s})\in V$
$$
(\phi_{(\cdot)}(S_{n}u_{s}),n(1-S_{n})v^{*}_{s})_{H}=
(\phi_{(\cdot)}(S_{n}u_{s}),S_{n}v^{*}_{s})_{V}.
$$
With probability one $u_{s}\in V$ for almost all $s$ 
for which also $\phi_{(\cdot)}(S_{n}u_{s})\to \phi_{(\cdot)}( u_{s})$
weakly in $V$, owing to assumption (iv), whereas
$S_{n}v^{*}_{s}\to v^{*}_{s}$ strongly in $V$. Hence with
probability one  for almost all $s$
$$
(\phi_{(\cdot)}(S_{n}u_{s}),n(1-S_{n})v^{*}_{s})_{H}\to
(\phi_{(\cdot)}( u_{s}), v^{*}_{s})_{V}.
$$
as $n\to\infty$. Furthermore,
$$
|(\phi_{(\cdot)}(S_{n}u_{s}),n(1-S_{n})v^{*}_{s})_{H}|=
|(\phi_{(\cdot)}(S_{n}u_{s}),S_{n}v^{*}_{s})_{V}| 
\leq K(1+\|u_{s}\|_{V})\|v^{*}\|_{V}
$$
by assumption (iii). It follows by the dominated convergence
theorem that (a.s.)
$$
\int_{0}^{t}(\phi_{(\cdot)}(S_{n}u_{s}),n(1-S_{n})v^{*}_{s})_{H}\,ds
\to \int_{0}^{t}(\phi_{(\cdot)}(u_{s}), v^{*}_{s})_{V}\,ds.
$$

Finally, (a.s.) for all $s$
$$
\phi_{(S_{n}\sigma^{k}_{s})(S_{n}\sigma^{k}_{s})}
(S_{n}u_{s})\to
\phi_{( \sigma^{k}_{s})( \sigma^{k}_{s})}
( u_{s})
$$
because of assumption (i) and
$$
|\phi_{(S_{n}\sigma^{k}_{s})(S_{n}\sigma^{k}_{s})}
(S_{n}u_{s})|\leq K(\sup_{s\leq t}\|u_{s}\|_{H})
\|\sigma^{k}_{s}\|_{H}^{2}
$$
in light of assumption (ii). This allows us to pass to the 
limit in the remaining expression in \eqref{8.6.5}
and brings the proof of the theorem to an end.

\mysection{The maximum principle for second-order SPDEs}
                                        \label{section 8.6.1}
In Section \ref{section 8.3.1} take a domain
$G\subset\bR^{d}$, $V=\HO^{1}_{2}=\HO^{1}_{2}(G)$, and 
$H=L_{2}=L_{2}(G)$. 

Take an  infinitely 
 differentiable function $r(x)$, $x\in\bR$, such 
  that $|r (x)|\leq N|x|^{2}$,
$|r' (x)|\leq N|x|$, and $|r'' |\leq N$, where $N$
is a constant. For $h\in L_{2}$ define
$$
\phi(h)=\int_{G}r(h(x))\,dx.
$$
As is easy to see,   assumptions (i) and (ii)
 of Section \ref{section 8.3.1}
are satisfied and for $h,\xi\in L_{2}$
$$
\phi_{(\xi)}(h)=\int_{G}r'(h(x))\xi(x)\,dx,\quad
\phi_{(\cdot)}(h)=r'(h(x)),
$$
$$
\phi_{(\xi)(\xi)}(h)=\int_{G}r''(h(x))\xi^{2}(x)\,dx.
$$
Furthermore, if $h\in\HO^{1}_{2}$, then there exists a sequence
of $h_{n}\in C^{\infty}_{0}$ such that
$h_{n}\to h$ in the norm of $\HO^{1}_{2}$. Almost obviously
$\phi_{(\cdot)}(h_{n})=r'(h_{n}(x))\in C^{\infty}_{0}$ 
and $r'(h_{n}(x))\to r'(h (x))$ in the $\HO^{1}_{2}$-norm.
Hence $\phi_{(\cdot)}(h )\in \HO^{1}_{2}$ if $h\in\HO^{1}_{2}$.
One can also easily verify that
$$
\|\phi_{(\cdot)}(h )\|_{H^{1}_{2}}
\leq N\|h\|_{H^{1}_{2}},
$$
where $N$ is the constant from above, so that
assumption   (iii)
 of Section \ref{section 8.3.1}
is satisfied as well. Finally, it is not hard to check that
for $v^{*}\in\HO^{1}_{2}$
$$
(r'(h ),v^{*} )_{H^{1}_{2}}=
\int_{G}r'(h(x))v^{*}(x)\,dx+\sum_{|\alpha|=1}
\int_{G}r''(h(x))(D_{i}h(x))D_{i}v^{*}(x)\,dx
$$
is continuous as a function of $h$ on the space $\HO^{1}_{2}$
and this is what is required in assumption (iv)
of  Section \ref{section 8.3.1}.

By Theorem \ref{theorem 8.3.1} we now conclude that
$$
 \int_{G}r(u_{t})\,dx=
\int_{G}r(u_{0})\,dx+\sum_{k}\int_{0}^{t}
\int_{G}r'(u_{s})\sigma^{k}_{s}\,dx\,dw^{k}_{s}
$$
$$
+\int_{0}^{t}\int_{G}r'(u_{s})v^{*}_{s}\,dx\,ds
+\sum_{|\alpha|=1}\int_{0}^{t}\int_{G}
r''(u_{s})(D_{i}u_{s}) D_{i}v^{*}_{s}\,dx\,ds
$$
\begin{equation}
                                             \label{8.6.6}
+(1/2)\int_{0}^{t}\int_{G}r''(u_{s})\sum_{k}|\sigma^{k}_{s}|^{2}
\,dx\,ds.
\end{equation}

Next, we generalize this formula for $r$ from a wider class.
Denote by $\cR$
 the set of real-valued  functions
$r(x)$ on $\bR$ such that

\noindent(i) $r$ is continuously differentiable, $r(0)=r'(0)=0$,

\noindent(ii) $r'$ is absolutely continuous, its
derivative $r''$ is bounded and left continuous,
that is usual $r''$ which exists almost everywhere
is bounded and there is a left-continuous
function with which $r''$ coincides
almost everywhere.

For $r\in\cR$ by $r''$ we will always mean the left-continuous
modification of the usual second-order derivative of $r$.

It turns out (see Remark 2.1 in \cite{Kr07}) that for any $r\in\cR$
 there exists a sequence
$r_{n}\in\cR$ of  infinitely 
 differentiable functions
 such that $|r_{n}(x)|\leq N|x|^{2}$,
$|r'_{n}(x)|\leq N|x|$, and $|r''_{n}|\leq N$ with $N
<\infty$ independent of
$x\in\bR$ and $n$, and
$r_{n},r_{n}',r_{n}''\to r,r'r''$ on $\bR$. By using this fact
one easily shows that $r(u)\in \HO^{1}_{2}$ if
$u\in \HO^{1}_{2}$ and \eqref{8.6.6} also holds for
$r\in\cR$.

In particular, we can apply \eqref{8.6.6}
with $r(x)=(x^{+})^{2}$ and by using the well-known fact that
\begin{equation}
                                               \label{8.7.8}
D_{i}(u^{+})=I_{u>0}D_{i}u
\end{equation}
we then obtain that
$$
\int_{G}(u^{+}_{t})^{2}\,dx
=\int_{G}(u^{+}_{0})^{2}\,dx +2\sum_{k}\int_{0}^{t}
\int_{G} u^{+}_{s}\sigma^{k}_{s}\,dx\,dw^{k}_{s}
$$
$$
+2\int_{0}^{t}\int_{G} u^{+}_{s} v^{*}_{s}\,dx\,ds
+2\sum_{i=1}^{d}\int_{0}^{t}\int_{G}
 (D_{i}u^{+}_{s}) D_{i}v^{*}_{s}\,dx\,ds 
$$
$$+
\int_{0}^{t}I_{u_{s}>0}\sum_{k}|\sigma^{k}|^{2}\,dx\,ds
=\int_{G}(u^{+}_{0})^{2}\,dx +2\sum_{k}\int_{0}^{t}
(u^{+}_{s},\sigma^{k}_{s})_{L_{2}}\,dw^{k}_{s}
$$
\begin{equation}
                                                   \label{8.7.4}
+2\int_{0}^{t} ( u^{+}_{s} ,v^{*}_{s})_{H^{1}_{2}}\,ds
+
\int_{0}^{t}I_{u_{s}>0}\sum_{k}|\sigma^{k}|^{2}\,dx\,ds.
\end{equation}

Next, assume that  in addition to \eqref{8.7.1} we have
that for any $\phi\in\HO^{1}_{2}$ (a.s.) for all $t$
$$
(\phi,u_{t })=(\phi,u_{0})+
\int_{0}^{t}
(\phi,\sigma^{ik}_{s}D_{i}u_{s}+\nu^{k}_{s}u_{s} )_{L_{2}}
 \,dw^{k}_{s}
$$
\begin{equation}
                                               \label{1.4.1}
+\int_{0}^{t}\big[(D_{i}\phi, -a^{ij}_{s}D_{j}u_{s}
-a ^{i}_{s}u_{s})_{L_{2}}
+(\phi,b^{i}_{s}D_{i}u_{s}+c_{s} u_{s} +f_{s})_{L_{2}}\big]
\,ds,
\end{equation}
where the summation with respect to repeated indices
is understood. We assume that $a^{ij}_{t}(x)$, $b^{i}_{t}(x)$, $a ^{i}_{t}(x)$,
$c_{t}(x)$, $f_{t}(x)$,
$\sigma^{ik}_{t}(x)$, and $\nu^{k}_{t}(x)$ 
 are real-valued  
functions defined for $i,j=1,...,d$, $k=1,2,...$,
$t\in[0,\infty)$, $x\in\bR^{d}$ and also depending on
$\omega\in\Omega$.

\begin{assumption}
                                          \label{assumption 1.10.1} 
For all values of the arguments

(i)   $\sigma^{i}:=(\sigma^{i1},\sigma^{i2},...)$,
$\nu:=(\nu^{1},\nu^{2},...) \in\ell_{2}$;

(ii) for all $\lambda\in\bR^{d}$ 
$$
 (2a^{ij}-\alpha^{ij})\lambda^{i}\lambda^{j}\geq0,
$$
 where
$\alpha^{ij}=(\sigma^{i},\sigma^{j})_{\ell_{2}}$.

\end{assumption}
Assumption \ref{assumption 1.10.1} (ii) is just the
usual parabolicity assumption.

We need one more function $K_{t}\geq0$ defined on $\Omega\times[0,
\infty)$.

\begin{assumption}
                                          \label{assumption 1.12.1}

(i) The functions
$a^{ij}_{t}(x)$, $b^{i}_{t}(x)$, $a ^{i}_{t}(x)$,
$c_{t}(x)$,
$\sigma^{ik}_{t}(x)$, $\nu^{k}_{t}(x)$,
  and $K (t)$ are measurable
with respect to
$(\omega,t,x)$ and $\cF_{t}$-adapted for each~$x$;

(ii)  the functions
$a^{ij}_{t}(x)$, $b^{i}_{t}(x)$, $a ^{i}_{t}(x)$,
$c_{t}(x)$,
$\sigma^{ik}_{t}(x)$, and $\nu^{k}_{t}(x)$ 
 are bounded;

(iii) for each $\omega,t$
the functions
$$
\eta^{i}_{t}:=a ^{i}_{t}-b^{i}_{t}
 -(\sigma^{i}_{t},\nu_{t})_{\ell_{2}}
$$
are once continuously differentiable on $D$,
have bounded derivatives, and satisfy
\begin{equation}
                                                \label{1.10.2}
D_{i}\eta^{i}+2c+|\nu|^{2}_{\ell_{2}}\leq K 
\end{equation}
for all values of arguments;

(iv)  
the process  $ f_{t}$, is   $L_{2} $-valued  
$\cF_{t}$-adapted and jointly measurable; and
 for all $T\in[0,\infty)$  
$$
E\int_{0}^{T}\big(  \|f_{s}\|^{2}_{L_{2} } 
+K_{s} \big)\,ds<\infty.
$$

\end{assumption}

 Under these assumptions (and the assumption that 
$u_{t}$ is taken from Section \ref{section 8.3.1}
corresponding to some $v^{*}_{t}$ and $m_{t}$) the stochastic
integrals in \eqref{1.4.1} have exactly the same form
as in \eqref{8.7.1} if in the latter we replace $\sigma^{k}_{s}$
with
\begin{equation}
                                                      \label{8.7.6}
\sum_{i=1}^{d}\sigma^{ik}_{s}D_{i}u_{s}+\nu^{k}_{s}u_{s} ,
\end{equation}
for which
\begin{equation}
                                                      \label{8.7.7}
\sum_{k} \big|\sum_{i=1}^{d}\sigma^{ik}_{s}D_{i}u_{s}+\nu^{k}_{s}u_{s}
\big |^{2} =
  \alpha^{ij}_{s}(D_{i}u_{s})D_{j}u_{s} 
+2 (\sigma^{i}_{s},\nu_{s})_{\ell_{2}}u_{s}D_{i}u_{s} 
+ |\nu_{s}|^{2}_{\ell_{2}}u^{2}_{s}.
\end{equation}
The processes \eqref{8.7.6}
  are predictable $L_{2}$-valued processes satisfying
$$
E\sum_{k}\int_{0}^{T}\|\sum_{i}
\sigma^{ik}_{s}D_{i}u_{s}+\nu^{k}_{s}u_{s} \|^{2}_{L_{2}}\,ds
$$
$$
\leq NE\sum_{i} \int_{0}^{T}\int_{G}
\alpha^{ii}_{s}(D_{i}u)^{2}_{s}\,dx \,ds
+NE\int_{0}^{T}\int_{G}
|\nu_{s}|^{2}_{\ell_{2}} u ^{2}_{s}\,dx \,ds<\infty
$$
 for any $T\in(0,\infty)$, where $N$ are  absolute constants.
 
At the first sight, the usual integral in \eqref{1.4.1}
does not look like the one in \eqref{8.7.1}. However, observe that
  on  $A:=\{(\omega,s):u_{s}\in
\HO^{1}_{2}\}$  the function 
$$
(D_{i}\phi, -a^{ij}_{s}D_{j}u_{s}
-a ^{i}_{s}u_{s})_{L_{2}}
+(\phi,b^{i}_{s}D_{i}u_{s} +c_{s}u_{s}+f_{s})_{L_{2}}
$$
as a function on $\HO^{1}_{2}$ is continuous. By  Riesz's
representation theorem there exists a unique $v^{*}$
such that 
\begin{equation}
                                                   \label{8.7.2}
I_{A}(D_{i}\phi, -a^{ij}_{s}D_{j}u_{s}
-a ^{i}_{s}u_{s})_{L_{2}}
+(\phi,b^{i}_{s}D_{i}u_{s}+c_{s} u_{s} +f_{s})_{L_{2}}=
(\phi,v^{*}_{s})_{H^{1}_{2}}.
\end{equation}
Since $u_{s}$ is an $L_{2}$-continuous process and
$\HO^{1}_{2}$ is a Borel subset of $L_{2}$, the set $A$
is predictable. Also notice that $D_{j}u_{s}$ could be defined
as the limits of finite differences. Hence, $I_{A}D_{j}u_{s}$
are also predictable and formula \eqref{8.7.2} (along
with the separability
of $\HO^{1}_{2}$) shows that $v^{*}_{s}$ is an $\HO^{1}_{2}$-valued
predictable process. Furthermore, the absolute value
of the left-hand side of \eqref{8.7.2} is obviously less than
$$
NI_{A}\|\phi\|_{H^{1}_{2}}(\|u_{s}\|_{H^{1}_{2}}
+\|f_{s}\|_{L_{2}}),
$$
where $N$ depends only on $d$ and the sup norms of the coefficients.
It follows that 
$$
E\int_{0}^{T}\|v^{*}_{s}\|^{2}_{H^{1}_{2}}\,ds<\infty
$$
for any $T$. 

Summing up all the above comments on equation 
\eqref{1.4.1} we conclude
that, our assumption that $u_{t}$ satisfies it,
 is justified if it satisfies \eqref{8.7.1} with $\sigma$ and $v^{*}$
specified above. We are not going to discuss the possibility
of existence of such $u_{t}$, that is the existence
of solutions of \eqref{1.4.1} in the class of functions $u_{t}$
as in Theorem \ref{theorem 4.9.3}. By the way, generally,
such solutions may not even exist. For instance, if
all the coefficients  and $f_{s}$ in \eqref{1.4.1}
vanish identically, we have $u_{t}=u_{0}$ and, if $u_{0}\notin
\HO^{1}_{2}$, we do not have $u_{t}\in \HO^{1}_{2}$
for almost all $(\omega,t)$.

We will just assume that we are given a
continuous $L_{2}$-valued predictable process
$u_{t}$ such that
$$
E\int_{0}^{T}\|u_{t}\|^{2}_{H^{1}{2}}\,dt<\infty
$$
for any $T\in[0,\infty)$ and equation \eqref{1.4.1}
holds for any $\phi\in\HO^{1}_{2}$ (a.s.) for all $t$.

A very particular case of the following theorem
can be found in \cite{Pa75}
(see also the references therein).

\begin{theorem}[maximum principle]
                                            \label{theorem 8.7.1}
Under the above assumptions suppose that $u_{0}\leq0$
and $f_{t}\leq0$ for almost all $(\omega,t)$.
Then (a.s.) for all $t$ we have $u_{t}\leq0$.

\end{theorem} 

Proof. According to what  has been explained before the theorem
and formulas \eqref{8.7.4}, \eqref{8.7.7},
 and \eqref{8.7.2} we have that
(a.s.) for all $t$
$$
\|u^{+}_{t}\|_{L_{2}}=
M_{t}+2\int_{0}^{t}
\big[(D_{i}u^{+}_{s}, -a^{ij}_{s}D_{j}u_{s}
-a ^{i}_{s}u_{s})_{L_{2}}
+(u^{+}_{s},b^{i}_{s}D_{i}u_{s}+c_{s} u_{s} +f_{s})_{L_{2}}\big]\,ds
$$
\begin{equation}
                                            \label{8.7.9}
+ \int_{0}^{t}\int_{G}I_{u_{t}>0}\big[
\alpha^{ij}_{s}(D_{i}u_{s})D_{j}u_{s} 
+2 (\sigma^{i}_{s},\nu_{s})_{\ell_{2}}u_{s}D_{i}u_{s} 
+ |\nu_{s}|^{2}_{\ell_{2}}u^{2}_{s}\big]\,dx\,ds
\end{equation}
where $M_{t}$ is a martingale.

According to \eqref{8.7.8}
$$
(a^{ij}_{s}(D_{i}u^{+}_{s}),D_{j}u_{s})_{L_{2}}=
( a^{ij}_{s}(D_{i}u^{+}_{s}),D_{j}u^{+}_{s})_{L_{2}}
$$
$$
(I_{u_{t}>0}\alpha^{ij}_{s}(D_{i}u_{s}),D_{j}u_{s})_{L_{2}}
=( \alpha^{ij}_{s}(D_{i}u^{+}_{s}),D_{j}u^{+}_{s})_{L_{2}},
$$
so that
$$
-2(a^{ij}_{s}(D_{i}u^{+}_{s}),D_{j}u_{s})_{L_{2}}
+(I_{u_{t}>0}\alpha^{ij}_{s}(D_{i}u_{s}),D_{j}u_{s})_{L_{2}}
\leq0
$$
in light of Assumption \ref{assumption 1.10.1} (ii).
Furthermore, at points where $u^{+}_{s}\in \HO^{1}_{2}$,  we have 
$$
I_{s}:=-(D_{i}u^{+}_{s},a^{i}_{s}u_{s})_{L_{2}}+
(u^{+}_{s},b^{i}_{s}D_{i}u_{s} )_{L_{2}}
+ ( (\sigma^{i}_{s},\nu_{s})_{\ell_{2}}u^{+}_{s},D_{i}u_{s})_{L_{2}}
$$
$$
=-(\eta^{i}_{s}u^{+}_{s},D_{i}u^{+}_{s})_{L_{2}}=(D_{i}\eta^{i}_{s}
u^{+}_{s},u^{+}_{s})_{L_{2}},
$$
where the last equality is obtained by integrating by parts,
which is justified by approximating $u^{+}_{s}\in \HO^{1}_{2}$
by $C^{\infty}_{0}$-functions and passing to the limit.
At this point the reader can understand that, actually,
we only need $\eta^{i}_{t}$ to be Lipschitz continuous rather than
continuously differentiable. In any case by also observing that
$u^{+}_{s}f_{s}\leq0$ and
using Assumption
\ref{assumption 1.12.1} (iii) we conclude from \eqref{8.7.9}
that (a.s.)
 $$
d\|u^{+}_{t}\|_{L_{2}}\leq
dM_{t}+K_{t}\|u^{+}_{t}\|_{L_{2}},
$$
$$
d\big[\|u^{+}_{t}\|_{L_{2}}\exp(-\int_{0}^{t}K_{s}\,ds)\big]
\leq\exp(-\int_{0}^{t}K_{s}\,ds)\,dM_{t},
$$
$$
\|u^{+}_{t}\|_{L_{2}}\exp(-\int_{0}^{t}K_{s}\,ds)
\leq\int_{0}^{t}
\exp(-\int_{0}^{s}K_{r}\,dr)\,dM_{r}.
$$
In the last relation the left-hand side is nonnegative and
the right-hand side is a martingale starting from  zero.
It follows that, with probability one,
 the martingale is zero and so is 
$\|u^{+}_{t}\|_{L_{2}}$ which proves the theorem.


\begin{thebibliography}{mm}

\bibitem{BG1} S. Bonaccorsi and G. Guatteri,
{\em Stochastic partial differential equations in bounded
domains with Dirichlet boundary conditions\/},
Stochastics and Stoc. Rep., Vol. 74 (1-2) (2002), 349-370.

\bibitem{BG2} S. Bonaccorsi and G. Guatteri,
{\em Classical solutions for SPDEs with 
Dirichlet boundary conditions\/}, pp. 33-44 in
Progress in Probability,
Vol. 52, Birkh\"auser, Basel/Switzerland, 2002.

\bibitem{PT}
G. Da Prato and L. Tubaro, {\em
 Fully nonlinear stochastic partial
differential equations\/},  SIAM J. Math. Anal., Vol. 27  
(1996),  No. 1, 40-55.



\bibitem{Br}
Z. Brze\'zniak, J. M. A. M.~van Neerven, M. C. Veraar, and L. Weis, {\em
Ito's formula in UMD Banach spaces and regularity of solutions of the Zakai
equation\/},  
J. Differential Eq., Vol. 245 (2008), 30-58.

\bibitem{Ev} L.C. Evans, ``Partial Differential Equations",
Graduate Studies in Mathematics, Vol. 19,
American Mathematical Society, Providence, RI, 1998.

\bibitem{GK} I. Gy\"ongy and N.V. Krylov,
  {\em  On stochastic equations with
respect to semimartingales II. It\^o formula in Banach spaces\/},
   Stochastics, Vol. 6   (1982), No. 3--4,   153--173.

\bibitem{Kr07} N.V. Krylov, {\em
Maximum principle for SPDEs 
and its applications\/}, pp. 311-338  in
``Stochastic Differential Equations: Theory and Applications,
A Volume in Honor of Professor Boris L. Rozovskii",
 P.H. Baxendale, S.V. Lototsky eds., Interdisciplinary 
Mathematical Sciences, Vol. 2, World Scientific, 2007.

\bibitem{KR81} N.V. Krylov and
 B.L. Rozovsky, {\em  On the first integrals and
Liouville equations for diffusion processes\/},
pp.  117-125 in ``Stochastic Differental Systems,
Proc. 3rd IFIP-WG 7/1 Working Conf., 
Visegr\'ad, Hungary, Sept.
15-20, 1980",
    Lecture Notes in Contr. Inform. Sci.,
   Vol. 36  (1981).

 

\bibitem{KR} N.V. Krylov and B.L. Rozovsky,  {\em  
Stochastic evolution
equations},  ``Itogy nauki i
tekhniki'',  Vol. 14, VINITI, Moscow, 1979, 71-146 in Russian;
English translation in
 J. Soviet Math., Vol. 16 (1981), No. 4, 1233-1277. 


\bibitem{KR82} N.V. Krylov and
 B.L. Rozovsky,   {\em
   On the characteristics of degenerate second order
parabolic It\^o equations}, 
 Trudy seminara imeni Petrovskogo,
  Vol. 8 (1982),  153-168 in Russian; English 
translation:
   J. Soviet Math Vol. 32
  (1986), No. 4,
   336-348.

\bibitem{KW} N.V. Krylov and Teng Wang, {\em Filtering partially
observable diffusions up to the exit time  from a domain\/}, 
Stoch. Proc. Appl., Vol. 121 (2011), No. 8, 1785--1815.

\bibitem{Ku} H. Kunita,
{\em On backward stochastic differential equations\/}, 
Stochastics, Vol.  6  (1981/82), No. 3-4, 293-313. 

\bibitem{Li} J.-L. Lions,
``Quelques methodes de resolution des problemes aux limites
non lineaire'', Dunod, Gauthier Villars, 1969.

\bibitem{Pa75} E.~Pardoux, {\em Equations aux d\'erivees
partielles  stochastiques non lin\'eaires monotones\/}, Ph.
D. Thesis, Universit\'e de Paris  Sud, Orsay, 1975,
http://www.cmi.univ-mrs.fr/~pardoux/Pardoux\_these.pdf

\bibitem{Pa79} E.~Pardoux, {\em Stochastic partial differential 
equations and filtering
of diffusion processes\/}, Stochastics, Vol. 3 (1979), No. 2,
 127-167.

\bibitem{Tu} L. Tubaro, {\em
Some results on stochastic partial differential
equations by the stochastic characteristics method\/},  Stochastic
Anal. Appl., Vol. 6  (1988),  No. 2, 217-230.

\end{thebibliography}
\end{document}